
\input amstex.tex
\documentstyle{amsppt}
\def\DJ{\leavevmode\setbox0=\hbox{D}\kern0pt\rlap
{\kern.04em\raise.188\ht0\hbox{-}}D}
\footline={\hss{\vbox to 2cm{\vfil\hbox{\rm\folio}}}\hss}
\nopagenumbers
\font\ff=cmr8

\baselineskip=13pt
\def\hf{{\textstyle{1\over2}}}
\def\a{\alpha}

\def\d{{\,\roman d}}
\def\e{\varepsilon}

\def\={\;=\;}

\def\z{\zeta}
\def\e{\varepsilon}

\font\teneufm=eufm10
\font\seveneufm=eufm7
\font\fiveeufm=eufm5
\newfam\eufmfam
\textfont\eufmfam=\teneufm
\scriptfont\eufmfam=\seveneufm
\scriptscriptfont\eufmfam=\fiveeufm
\def\mathfrak#1{{\fam\eufmfam\relax#1}}

\font\tenmsb=msbm10
\font\sevenmsb=msbm7
\font\fivemsb=msbm5
\newfam\msbfam
\textfont\msbfam=\tenmsb
\scriptfont\msbfam=\sevenmsb
\scriptscriptfont\msbfam=\fivemsb
\def\Bbb#1{{\fam\msbfam #1}}

\def \NN {\Bbb N}

\def\rightheadline{{\hfil{\ff
Sum involving $\pi(x)$}\hfil\tenrm\folio}}

\def\leftheadline{{\tenrm\folio\hfil{\ff
Aleksandar Ivi\'c }\hfil}}
\def\emptyheadline{\hfil}
\headline{\ifnum\pageno=1 \emptyheadline\else
\ifodd\pageno \rightheadline \else \leftheadline\fi\fi}

\topmatter
\title ON A SUM INVOLVING THE PRIME COUNTING FUNCTION
$\pi(x)$ \endtitle

\author   Aleksandar Ivi\'c
\endauthor
\dedicatory Univ. Beograd. Publikac. Elektrotehn. Fak. Ser. Mat.
{\bf13}${\textstyle {(2002), 85-88}}$
\enddedicatory
\address{
Aleksandar Ivi\'c, Katedra Matematike RGF-a
Universiteta u Beogradu, \DJ u\v sina 7, 11000 Beograd,
Serbia (Yugoslavia).}
\endaddress
\email {\tt eivica\@ubbg.etf.bg.ac.yu,
aivic\@rgf.bg.ac.yu} \endemail

\abstract
{An asymptotic formula for the sum of reciprocals of $\pi(n)$
is derived, where $\pi(x)$ is the number of primes not exceeding $x$.
This result improves the previous results of De Koninck-Ivi\'c and
L. Panaitopol.}
\endabstract

\endtopmatter
Let, as usual, $\pi(x) = \sum_{p\leqslant x}1$ denote the
number of primes not exceeding $x$. The prime number theorem
(see e.g., [2, Chapter 12]) in its strongest known form states that
$$
\pi(x) \= \text{li}\,x + R(x),\leqno(1)
$$
with
$$
\text{li}\,x := \int_2^x{\d t\over\log t} =
x\left({1\over\log x} +   {1!\over\log^2 x}  + \cdots+
{m!\over\log^{m+1} x} +  O\left({1\over\log^{m+2} x}\right)\right)\leqno(2)
$$
for any fixed integer $m \geqslant 0$, and
$$
R(x) \ll x\exp(-C\delta(x)),\quad \delta(x) := (\log x)^{3/5}
(\log\log x)^{-1/5}\quad(C>0),\leqno(3)
$$
where henceforth $C,C_1,\ldots\,$ will denote absolute constants.
In [1, Theorem 9.1] J.-M. De Koninck and the author proved that
$$
\sum_{2\leqslant n\leqslant x}{1\over\pi(n)} = \hf\log^2x + O(\log x).
\leqno(4)
$$
Recently L. Panaitopol [1] improved (4) to
$$
\sum_{2\leqslant n\leqslant x}{1\over\pi(n)} = \hf\log^2x - \log x
-\log\log x + O(1).
\leqno(5)
$$
One obtains (5) from the asymptotic formula
$$
{1\over\text{li}\,x} = {1\over x}\left(\log x - 1 - {k_1\over\log x}
- {k_2\over\log^2 x} - \ldots - {k_m(1+\a_m(x))\over\log^m x}\right),
\leqno(6)
$$
where $\,\a_m(x) \ll_m 1/\log x\,$, and the constants $\,k_1,\ldots\,,k_m\,$
are defined by the  recurrence  relation
$$
k_m +1!k_{m-1} + \ldots + (m-1)!k_1 =  m\cdot m!\qquad(m \in \NN),\leqno(7)
$$
so that $k_1 = 1,\,k_2 = 3,\, k_3 = 13,$ etc. This was established in [3].
Using (6) we shall give a further improvement of (4),
contained in the following

\bigskip
THEOREM. {\it For any fixed integer $m \geqslant 2$ we have}
$$
\eqalign{
\sum_{2\leqslant n\leqslant x}{1\over\pi(n)} &= \hf\log^2x - \log x
-\log\log x + C\cr&
+ {k_2\over\log x} +{k_3\over2\log^2 x} + \ldots +
{k_m\over(m-1)\log^{m-1} x} + O\left({1\over\log^m x}\right),
\cr}\leqno(8)
$$
{\it where $C$ is an absolute constant, and $k_2,\,\ldots\,,k_m$
are the constants defined by }(7).

\bigskip
{\bf Proof}. From (1) we have
$$
\eqalign{
\sum_{2\leqslant n\leqslant x}{1\over\pi(n)} &=
1 + \sum_{3\leqslant n\leqslant x}{1\over\text{li}\,n}
- \sum_{3\leqslant n\leqslant x}{R(n)\over\text{li}\,n
(\text{li}\,n + R(n))}\cr&
= \sum_{3\leqslant n\leqslant x}{1\over\text{li}\,n}
+ \left\{1 - \sum_{n=3}^\infty {R(n)\over\text{li}\,n
(\text{li}\,n + R(n))}\right\} +
\sum_{n> x}{R(n)\over\text{li}\,n
(\text{li}\,n + R(n))}\cr&
= \sum\nolimits_1\, + \,C_1\, + \,\sum\nolimits_2,\cr}
$$
say. By using the bound  $\,\text{li}\,x \ll x/\log x\,$ and (3) it is seen that
$$
\eqalign{
\sum\nolimits_2 &= \sum_{n> x}{R(n)\over\text{li}\,n
(\text{li}\,n + R(n))} \ll \sum_{n>x}{1\over n}
{\roman e}^{-{1\over2}C\delta(n)}\cr&
\ll {\roman e}^{-{1\over3}C\delta(x)}\int_{x-1}^\infty{1\over t}
{\roman e}^{-{1\over6}C\delta(t)}\d t \ll
{\roman e}^{-{1\over3}C\delta(x)},\cr}
$$
since $\delta(x)$ is increasing for $x\geqslant x_0$, and the
substitution $\log t= u$ easily shows that the above integral
is convergent. To evaluate $\sum_1$ we need the familiar Euler--Maclaurin
summation formula (see e.g., [2, eq. (A.23)]) in the form
$$
\sum_{X<n\leqslant X}f(n) = \int_X^Y f(t)\d t - \psi(Y)f(Y) + \psi(X)f(X)
+ \int_X^Y \psi(t)f'(t)\d t,\leqno(9)
$$
where $\psi(x) = x - [x] - \hf$ and $f(x) \in C^1[X,\,Y]\,$. We obtain from
(6), for $m \geqslant 2$ a fixed integer,
$$
\eqalign{
\sum\nolimits_1 &= \sum_{3\leqslant n\leqslant x}{1\over\text{li}\,n}\cr&
= \sum_{3\leqslant n\leqslant x}{1\over n}
\left(\log n - 1 - {k_1\over\log n} - {k_2\over\log^2 n} - \ldots
- {k_m(1+\a_m(n))\over\log^m n}\right),\cr}\leqno(10)
$$
and we evaluate each sum in (10) by using (9). We obtain
$$\eqalign{
\sum_{3\leqslant n\leqslant x}{\log n\over n} &= \hf\log^2x + c_1 +
O\left({\log x\over x}\right),\cr
\sum_{3\leqslant n\leqslant x}{1\over n} &= \log x + c_2 +
O\left({1\over x}\right),\cr
\sum_{3\leqslant n\leqslant x}{k_1\over n\log n} &= \log\log x + c_3
+ O\left({1\over x\log x}\right),\cr}
$$
and for $2 \leqslant r \leqslant m$
$$
\eqalign{
\sum_{3\leqslant n\leqslant x}{k_r\over n\log^r n} &=
k_r\int_3^x{\d t\over t\log^rt} + C_r + O\left({1\over x\log^r x}\right) \cr&
= k_r\int_3^\infty{\d t\over t\log^rt} -
k_r\int_x^\infty{\d t\over t\log^rt}
+ C_r + O\left({1\over x\log^r x}\right) \cr&
= D_r - {k_r\over(r-1)\log^{r-1}x} + O\left({1\over x\log^r x}\right) \cr}
$$
with
$$
D_r = C_r + k_r\int_3^\infty{\d t\over t\log^rt},
$$
Finally in view of $\,\a_m(x) \ll 1/\log x\,$ it follows that, for
$m\geqslant 2$ fixed,
$$
\sum_{3\leqslant n\leqslant x}{k_m\alpha_m(n)\over n\log^m n}
= \sum_{n=3}^\infty{k_m\alpha_m(n)\over n\log^m n} +
O\left({1\over\log^mx}\right).
$$
Putting together the above expressions in (10) we infer that
$$
\eqalign{
\sum\nolimits_1 &= \hf\log^2x - \log x
-\log\log x + C\cr&
+ {k_2\over\log x} +{k_3\over2\log^2 x} + \ldots +
{k_m\over(m-1)\log^{m-1} x} + O\left({1\over\log^m x}\right),
\cr}
$$
and then (8) easily follows with
$$
C = C_1 + c_1 - c_2 - c_3 - D_2 - \ldots - D_m - \sum_{n=3}^\infty\,
{k_m\a_m(n)\over n\log^mn}.
$$
The constant  $C$ in (8) does not
depend on $m$, which can be easily seen by taking two different
values of $m$ and then comparing the results.

\medskip
Note that we can evaluate directly $\sum_1$ by the Euler-Maclaurin
summation formula to obtain
$$
\sum\nolimits_1 = \int_3^x{\d t\over\text{li}\,t}
+ C_0 + O\left({\log x\over x}\right).\leqno(11)
$$
Integration by parts gives, for $x > 3$,
$$
\int_3^x{\d t\over\text{li}\,t} = \int_3^x\log t\d(\log\text{li}\,t) =
\log x\log(\text{li}\,x)
- \int_3^x{\log(\text{li}\,t)\over t}\d t - \log3\log\text{li}\,3,
$$
which inserted into (11) gives another expression for our sum, namely
$$
\sum_{2\leqslant n\leqslant x}{1\over\pi(n)}
= \log x\log(\text{li}\,x)
- \int_3^x{\log(\text{li}\,t)\over t}\d t + B
+ O\left({\roman e}^{-D\delta(x)}\right)\;(D > 0),\leqno(12)
$$
from which we can again deduce (8) by using (2). The advantage of (12)
is that it has a sharper error term than (8), but on the other hand
the  expressions on the right-hand side of (12) involve the
non-elementary function $\text{li}\,x$. Note also that the Riemann
hypothesis (that all complex zeros of the Riemann zeta-function $\z(s)$
have real parts equal to $\hf$) is equivalent to the statement (see [2])
that, for any given $\e > 0$, $R(x) \ll x^{1/2+\e}$ in (3),
which would correspondingly
improve the error term in (12) to $O(x^{-1/2+\e})$.

\bigskip

\Refs

\bigskip

\item{[1]} J.-M. De Koninck and A. Ivi\'c, {\it Topics
in Arithmetical Functions}, Notas de Matem\'atica {\bf43},
North-Holland, Amsterdam etc., 1980.

\item{[2]} A. Ivi\'c, {\it The Riemann zeta-function}, John
Wiley \& Sons, New York, 1985.

\item{[3]} L. Panaitopol, {\it A formula for $\pi(x)$ applied to a
result of Koninck-Ivi\'c}, Nieuw \ Archief \ voor\ Wiskunde {\bf5/1}(2000),
55-56.

\vskip2cm

Aleksandar Ivi\'c

Katedra Matematike RGF-a

Universitet u Beogradu

\DJ u\v sina 7, 11000 Beograd

Serbia and Montenegro

{\tt aivic\@matf.bg.ac.yu, aivic\@rgf.bg.ac.yu}

\endRefs

\bye